\documentclass[fleqn,reqno,11pt,a4paper,final]{amsart}

\usepackage[a4paper,left=30mm,right=30mm,top=30mm,bottom=30mm,marginpar=20mm]{geometry} 
\usepackage{amsmath}
\usepackage{amssymb}
\usepackage{amsthm}
\usepackage{amscd}
\usepackage[ansinew]{inputenc}
\usepackage{cite}
\usepackage{bbm}
\usepackage{color}
\usepackage[english=american]{csquotes}
\usepackage[final]{graphicx}
\usepackage{hyperref}
\usepackage{calc}
\usepackage{mathptmx}

\graphicspath{{../Pictures/}}

\numberwithin{equation}{section}

\newtheoremstyle{thmlemcorr}{10pt}{10pt}{\itshape}{}{\bfseries}{.}{10pt}{{\thmname{#1}\thmnumber{ #2}\thmnote{ (#3)}}}
\newtheoremstyle{thmlemcorr*}{10pt}{10pt}{\itshape}{}{\bfseries}{.}\newline{{\thmname{#1}\thmnumber{ #2}\thmnote{ (#3)}}}
\newtheoremstyle{remexample}{10pt}{10pt}{}{}{\bfseries}{.}{10pt}{{\thmname{#1}\thmnumber{ #2}\thmnote{ (#3)}}}
\newtheoremstyle{ass}{10pt}{10pt}{}{}{\bfseries}{.}{10pt}{{\thmname{#1}\thmnumber{ A#2}\thmnote{ (#3)}}}

\theoremstyle{thmlemcorr}
\newtheorem{theorem}{Theorem}
\numberwithin{theorem}{section}

\newtheorem{proposition}[theorem]{Proposition}

\newtheorem{definition}[theorem]{Definition}

\theoremstyle{thmlemcorr*}
\newtheorem{theorem*}{Theorem}
\newtheorem{lemma*}[theorem]{Lemma}
\newtheorem{corollary*}[theorem]{Corollary}
\newtheorem{proposition*}[theorem]{Proposition}
\newtheorem{problem*}[theorem]{Problem}
\newtheorem{conjecture*}[theorem]{Conjecture}
\newtheorem{definition*}[theorem]{Definition}

\theoremstyle{remexample}
\newtheorem{remark}[theorem]{Remark}

\theoremstyle{ass}

\newcommand{\Rbb}{\mathbb{R}}
\newcommand{\Sbb}{\mathbb{S}}

\DeclareMathOperator{\supp}{supp}

\newcommand{\n}{\tilde{n}}

\def\XXint#1#2#3{{\setbox0=\hbox{$#1{#2#3}{\int}$} 
\vcenter{\hbox{$#2#3$}}\kern-.5\wd0}}

\renewcommand{\epsilon}{\varepsilon}
\renewcommand{\phi}{\varphi}

\begin{document}

\title[Renewal equation with measure data]{Generalized Entropy Method for the Renewal Equation with Measure Data}

\author{Piotr Gwiazda}
\address{\textit{Piotr Gwiazda:} Institute of Mathematics, Polish Academy of Sciences, \'Sniadeckich 8, 00-656 Warszawa, Poland, and Institute of Applied Mathematics and Mechanics, University of Warsaw, Banacha 2, 02-097 Warszawa, Poland}
\email{pgwiazda@mimuw.edu.pl}

\author{Emil Wiedemann}
\address{\textit{Emil Wiedemann:} Hausdorff Center for Mathematics and Mathematical Institute, Universit\"{a}t Bonn, Endenicher Allee 60, 53115 Bonn, Germany}
\email{emil.wiedemann@hcm.uni-bonn.de}

\begin{abstract}
We study the long-time asymptotics for the so-called McKendrick-Von Foerster or renewal equation, a simple model frequently considered in structured population dynamics. In contrast to previous works, we can admit a bounded measure as initial data. To this end, we apply techniques from the calculus of variations that have not been employed previously in this context. We demonstrate how the generalized relative entropy method can be refined in the Radon measure framework. 
\end{abstract}

\maketitle

{\bf Keywords:} Structured population model, positive Radon measures, generalized relative entropy methods, measure valued-solutions, concentration measure.\\
\section{Introduction} 
This paper is devoted to the study of the long-time asymptotics of a linear structured population model, the so-called McKendrick-Von Foerster equation (or renewal equation, in probabilistic description). The equation gives the simplest well-known structured population model, which we choose here in order to illustrate the usefulness of certain variational tools in proving the long-time asymptotics with measure initial data.   

Classical studies on this topic were restricted to initial data in $L^1$ and were based on the analysis of semigroup theory for positive, irreducible operators (see e.g.\ \cite{ArinoWebb} and the monograph \cite{Webb}), or followed Feller's approach for Markov Processes, making use of the Laplace transform (see e.g.\ the monograph \cite{Iannelli}. Both approaches give rise to an exponential convergence result under the assumption of a spectral gap property.
However, for the model discussed in this paper and in contrast to the selection-mutation equation \cite{BuergerBook, Jabin}, the presence of the transport term destroys strong continuity of the semigroup in the space of measures with total variation norm, so that new ideas are needed to include measure initial data.

An entirely new approach was proposed by B.\ Perthame and collaborators in the papers \cite{MischlerPerthameRyzhik,MichelMischlerPerthame,MichelMischlerPerthame2}, see also the subsequent monograph \cite{Perthame}.  
This method, called generalized relative entropy method, was based on multiplying the linear equation by some nonlinear function of the solution in order to obtain a family of nonlinear renormalizations (so-called relative entropies). 
A clear adventage of this method was to enable the proof of a convergence result even in the absence of a spectral gap.
However this method seems prima facie also restricted to solutions in $L^1$, since the composition of nonlinear function with a Radon measure has no obvious meaning.

On the other hand, in the last years there was a strong development of existence and uniqueness theory and the convergence of numerical schemes in spaces of nonegative measures for structure population models \cite{EBT,CanizoCC,CCGU,CGU,ColomboGuerra2009,DiekmannGetto,GLMC,GMC,GJMCU,two_sex}, as well as for related crowded dynamic models, see e.g.\ \cite{EversHilleMuntean,EversHilleMuntean2,piccoli}. The motivation for considering such measure solutions comes from the desire to treat discrete and continuous initial distributions in a unified way. The aim of this paper is thus to fill the gap in the existing theory of measure solutions in terms of the long-time asymptotics, and to extend the method of generalized relative entropy to the case of initial data in the space of measures.

The method presented below is influenced by recent studies on measure-valued-strong uniqueness (based on the relative energy method). This direction of research started with the incompressible Euler equations \cite{BrDLSz}, and was later extended to nonlinear elasticity \cite{DeStTz} as well as compressible Euler \cite{GwSwWi} and Navier-Stokes equations \cite{FGSGW}.

To deal with measure data, we use the notion of recession function, which allows in a sense to take a nonlinear function of a bounded measure, and a continuity theorem for certain functionals due to Reshetnyak \cite{reshetnyak} and later refined by Kristensen-Rindler \cite{KR}. It turns out (Theorem \ref{expasympthm}) that the combination of these techniques, which have not previously been exploited in population models, and known results for the $L^1$ setting allows for a remarkably simple proof of the long-time asymptotics for measure initial data.

In Section \ref{GRE}, we formulate and prove the generalized relative entropy inequality in the context of measure solutions and show how this approach yields an alternative way to show the long-time asymptotics. Of course, Theorem \ref{alternative} is weaker than Theorem \ref{expasympthm} and is proved in a more complicated fashion, so that Section \ref{GRE} seems redundant if we focus only on the McKendrick-Von Foerster model. However we expect that the generalized relative entropy approach will be fruitful also for more sophisticated structured population models, for which statements like Theorem \ref{expasympthm} (related to the hypercontractivity property and the spectral gap) are not available even for $L^1$ data. We hope to achieve such results in the future, but for the time being our aim is to demonstrate how, in the simplest case, the relative entropy method can be used in the measure setting.  

\section{The Model}
Following the presentation and notation from \cite{Perthame}, we consider the McKendrick-Von Foerster equation or renewal equation in the form 
\begin{equation}\label{originalrenewal}
\begin{aligned}
\partial_t n(t,x)+\partial_x n(t,x)&=0\hspace{0.3cm}\text{on $(\Rbb^+)^2$},\\
n(t,x=0)&=\int_0^\infty B(y) n(t,y)dy,\\
n(t=0,x)&=n^0(x).
\end{aligned}
\end{equation}
Here, $n(x,t)$ denotes the population density at time $t$ with age $x$, and $B\in L^\infty(\Rbb^+;\Rbb^+)$ is a birth rate with the property that there exists a $\lambda_0>0$ such that
\begin{equation*}
\int_0^\infty B(x)dx>1.
\end{equation*}
Under these assumptions, it can be shown that there exist uniquely determined solutions of the primal and dual eigenvalue problems,
\begin{equation}\label{primal}
\begin{aligned}
\partial_x N(x)+\lambda_0 N(x)&=0,\hspace{0.2cm}x\geq0,\\
N(0)&=\int_0^\infty B(y)N(y)dy,\\
N>0,&\hspace{0.3cm}\int_0^\infty N(x)dx=1
\end{aligned}
\end{equation}
and
\begin{equation*}
\begin{aligned}
-\partial_x\phi(x)+\lambda_0\phi(x)=\phi(0)B(x),\hspace{0.2cm}x\geq0,\\
\phi\geq0,&\hspace{0.3cm}\int_0^\infty N(x)\phi(x)dx=1.
\end{aligned}
\end{equation*}
In fact, one easily discovers that the solution of \eqref{primal} is given by $N(x)=\lambda_0 e^{-\lambda_0x}$. As the birth rate has integral greater than one and we do not include a death rate, one expects the population to grow exponentially in time. In order to quotient out this growth, we set 
\begin{equation*} 
\n(t,x)=n(t,x)e^{-\lambda_0t},
\end{equation*}
whereupon \eqref{originalrenewal} becomes

\begin{equation}\label{renewal}
\begin{aligned}
\partial_t\n(t,x)+\partial_x\n(t,x)+\lambda_0\n(t,x)&=0\hspace{0.3cm}\text{on $(\Rbb^+)^2$},\\
\n(t,x=0)&=\int_0^\infty B(y)\n(t,y)dy,\\
\n(t=0,x)&=n^0(x).
\end{aligned}
\end{equation}
The formulation presented above is obviously not valid in the case of initial data in the space of measures, because for weak (or distributional) solutions, a pointwise equation involving derivatives of the unknown function has no meaning. 
Therefore as usual solutions (which are then proved to be Lipschitz continuous in the space of measures equipped with a Lipschitz-bounded distance) are understood in an integral sense only; see for more details \cite{GLMC,GMC} (see e.g.\ Definition 3.1. in \cite{GLMC}). We also refer to these papers for a result on the Lipschitz dependence of the solution on the initial data (in the space of measures equipped with the Lipschitz-bounded distance, see e.g. proposition 3.9 \cite{GLMC}).
Since the bounded Lipschitz distance metrizes the weak-star convergence on balls with respect to the total variation norm (see e.g.\ Theorem 2.7 in \cite{GLMC}), this also implies that if we consider a sequence of initial data converging weakly-star, then also the corresponding sequence of solutions will converge weakly-star (for every fixed time).

\section{Recession Functions and Continuity of Functionals}
We are interested in the long time asymptotics as $t\to\infty$ of this equation in the case that the initial data $n^0$ is only a bounded measure on $[0,\infty)$. To study this problem we need some tools from the calculus of variations, which we recall in the sequel.

First, suppose $f:\Rbb^n\to\Rbb$ is a continuous function with at most linear growth: $|f(z)|\leq C(1+|z|)$. We define (if it exists) its \emph{recession function} as
\begin{equation*}
f^\infty(z)=\lim_{s\to\infty}\frac{f(sz)}{s},\hspace{0.3cm}z\in\Rbb^n\setminus\{0\}.
\end{equation*}
Note that $f^\infty$ is 1-homogeneous, i.e. $f^\infty(\alpha z)=\alpha f^\infty(z)$ for any $\alpha>0$, so that it is completely determined by its values on the unit sphere $\Sbb^{n-1}$.

\begin{definition}
The set $\mathcal{F}(\Rbb^n)$ of continuous functions $f:\Rbb^n\to\Rbb$ which have a recession function that is continuous on $\Sbb^{n-1}$ is called the class of \emph{admissible integrands}.
\end{definition}

Given a domain $\Omega\subset\Rbb^m$ and a (possibly vector-valued) finite measure $\gamma\in \mathcal{M}(\overline{\Omega};\Rbb^n)$, we can write its Radon-Nikod\'ym decomposition w.r.t.\ Lebesgue measure as $\gamma=\gamma^a(x)dx+\gamma^s$, where $\gamma^s$ and $dx$ are mutually singular. We write
\begin{equation*}
\langle\gamma\rangle:=\int_{\Omega}\sqrt{1+|\gamma^a|^2}dx+|\gamma^s|(\overline{\Omega}).
\end{equation*}
We have the following (semi-)continuity properties (the first one is well-known, the second one is known as \emph{Reshetnyak's continuity theorem} \cite{reshetnyak, KR}):
\begin{proposition}\label{continuity}
Let $\{\gamma_n\}$ be a bounded sequence in $\mathcal{M}(\overline{\Omega};\Rbb^n)$ and assume $\gamma_n\stackrel{*}{\rightharpoonup}\gamma$ weakly* in the space of measures. Let $f\in\mathcal{F}(\Rbb^n)$.
\begin{itemize}
\item[a)] If $f$ is convex and $\psi\in C_b(\overline{\Omega})$, $\psi\geq0$, then
\begin{equation*}
\begin{aligned}
\liminf_{n\to\infty}\left\{\int_\Omega\psi(x)f(\gamma_n^a(x))dx+\int_{\overline{\Omega}}\psi(x)f^\infty\left(\frac{\gamma_n^s}{|\gamma_n^s|}\right)d|\gamma_n^s|(x)\right\}&\\
\geq\int_\Omega\psi(x)f(\gamma^a(x))dx+&\int_{\overline{\Omega}}\psi(x)f^\infty\left(\frac{\gamma^s}{|\gamma^s|}\right)d|\gamma^s|(x).
\end{aligned}
\end{equation*}
\item[b)] If $\langle\gamma_n\rangle\to\langle\gamma\rangle$ and $\psi\in C_b(\overline{\Omega})$, then
\begin{equation*}
\begin{aligned}
\lim_{n\to\infty}\left\{\int_\Omega\psi(x)f(\gamma_n^a(x))dx+\int_{\overline{\Omega}}\psi(x)f^\infty\left(\frac{\gamma_n^s}{|\gamma_n^s|}\right)d|\gamma_n^s|(x)\right\}&\\
=\int_\Omega\psi(x)f(\gamma^a(x))dx+&\int_{\overline{\Omega}}\psi(x)f^\infty\left(\frac{\gamma^s}{|\gamma^s|}\right)d|\gamma^s|(x).
\end{aligned}
\end{equation*}
\end{itemize} 
\end{proposition} 

\section{Long Time Asymptotics}
It was proved in \cite{GLMC} that \eqref{renewal} has a unique solution in the sense of distributions when $n^0\in\mathcal{M}^+([0,\infty))$. We show the following result on the long-time behavior of this solution:
\begin{theorem}\label{expasympthm}
Let $n^0\in\mathcal{M}^+([0;\infty))$. Then there is $y_0>0$, $\sigma>0$ and a bounded function $\eta$, positive on $\supp\phi$, such that the solution of the renewal equation satisfies
\begin{equation}\label{expasymp}
\int_0^\infty \eta(x)d|\n(t,x)-m_0N(x)dx|\leq e^{-\sigma(t-y_0)} \int_0^\infty \eta(x)d|\n^0(x)-m_0N(x)dx|,
\end{equation}
where $m_0=\int_0^\infty\phi(x)dn^0(x)$.
\end{theorem}
\begin{proof}
We use the variational techniques from the previous section in order to argue by approximation. Let $n^0_\epsilon$ be a regularization of $n^0$ such that $n^0_\epsilon\stackrel{*}{\rightharpoonup}n^0$ in the sense of measures and $\langle n^0_\epsilon \rangle\to \langle n^0\rangle$. By Theorem 1.1 in \cite{GwiazdaPerthame},~\eqref{expasymp} holds true for $\n_\epsilon(t,x)$ (the solution emanating from $n^0_\epsilon$), with $m_0$ replaced by $m_\epsilon:=\int_0^\infty\phi(x)dn^0_\epsilon(x)$. Note carefully that $y_0$, $\sigma$, and $\eta$ do not depend on $\epsilon$. Moreover, for every $t>0$, $\n_\epsilon(t)\stackrel{*}{\rightharpoonup}\n(t)$ in the sense of measures; indeed this follows from Theorems 2.7 and 4.6 in~\cite{GLMC}. This immediately implies $m_\epsilon\to m_0$.

The right hand side of~\eqref{expasymp} with $\epsilon$ converges, as $\epsilon\to0$, to the right hand side with $n^0$ and $m_0$; indeed this follows from Proposition~\ref{continuity}b) setting $\gamma_\epsilon=\n(t,x)_\epsilon-m_\epsilon N(x)dx$, $f=|\cdot|$, and $\psi=\eta$. Likewise, by Proposition~\ref{continuity}a), the $\liminf$ of the left hand side, as $\epsilon\to0$, is no less than the left hand side with $\n$ and $m_0$. Thus, \eqref{expasymp} is already proved.  
\end{proof}

\begin{remark}
In fact the proof of Theorem \ref{expasympthm} is not specific to the renewal equation. Indeed, whenever a contraction property like \eqref{expasymp} is known for some model for $L^1$ data (even for a speed of convergence other than exponential), our approximation argument works, provided the existence of measure-valued solutions is available.
\end{remark}

\section{Generalized Relative Entropy}\label{GRE}

Similar techniques can be applied in order to formulate and prove the generalized relative entropy inequality.
\begin{theorem}
Let $\n(t,x)$ be the solution of~\eqref{renewal} with $n^0\in\mathcal{M}([0;\infty))$. 
\begin{itemize}
\item[a)] If $H:\Rbb\to\Rbb^+$ is a convex admissible integrand, then
\begin{equation}\label{contraction}
\frac{d}{dt}\left\{\int_0^\infty\phi(x)N(x)H\left(\frac{\n^a(t,x)}{N(x)}\right)dx+\int_0^\infty\phi(x)H^\infty\left(\frac{\n^s(t)}{|\n^s(t)|}(x)\right)|\n^s(t,dx)| \right\}\leq 0
\end{equation} 
in the sense of distributions.
\item[b)] Let $\mu=\frac{B(x)}{N(0)}N(x)dx$. If $H$ is a convex admissible integrand, then
\begin{equation}\label{jensenquant}
\begin{aligned}
\int_0^\infty\left\{\int_0^\infty H\left(\frac{\n^a(t,x)}{N(x)}\right)d\mu(x)+\int_0^\infty\frac{B(x)}{N(0)}H^\infty\left(\frac{\n^s(t)}{|\n^s(t)|}(x)\right)|\n^s(t,dx)|-H\left(\int_0^\infty\frac{B(x)}{N(0)}\n(t,dx)\right)\right\}dt&\\
\leq \int_0^\infty\phi(x)N(x)H\left(\frac{(n^0)^a(x)}{N(x)}\right)dx+\int_0^\infty\phi(x)H^\infty\left(\frac{(n^0)^s}{|(n^0)^s|}(x)\right)|(n^0)^s(dx)|.\hspace{1cm}&
\end{aligned}
\end{equation}

\end{itemize}
\end{theorem}
\begin{proof}
Let again $n^0_\epsilon$ be a regularization of $n^0$ with $n^0_\epsilon\stackrel{*}{\rightharpoonup}n^0$ in the sense of measures and $\langle n^0_\epsilon \rangle\to \langle n^0\rangle$. Theorem 3.3 in \cite{Perthame} and its proof imply that~\eqref{contraction} and~\eqref{jensenquant} are true for $\n_\epsilon(t,x)$, the solution arising from $n^0_\epsilon$. And again, for every $t>0$, $\n_\epsilon(t)\stackrel{*}{\rightharpoonup}\n(t)$ in the sense of measures.

Therefore, using the lower semicontinuity from Proposition~\ref{continuity}, the left hand side of~\eqref{contraction} is not greater than the $\liminf$ of the corresponding expressions for $\n_\epsilon$. On the other hand, the right hand side of~\eqref{contraction} for $\n_\epsilon$ converges to the right hand side for $\n$ by part b) of Proposition~\ref{continuity}. Altogether we obtain a).

Similarly, for~\eqref{jensenquant} the right hand side converges as $\epsilon\to0$ and the first two space integrals of the left hand side, for fixed $t>0$, might only decrease in the $\liminf$. Moreover, for the third integral,
\begin{equation*}
\int_0^\infty\frac{B(x)}{N(0)}\n_\epsilon(t,dx)\to\int_0^\infty\frac{B(x)}{N(0)}\n(t,dx)
\end{equation*}   
for every $t$ by the weak convergence of $\n_\epsilon(t)$.

For simplicity let us denote by $J(t)$ the time integrand in~\eqref{jensenquant} (so that $\int J(t)dt$ is the left hand side of \eqref{jensenquant}) and by $J_\epsilon(t)$ the corresponding expression with $\n$ replaced by $\n_\epsilon$. We have shown that $\liminf_{\epsilon\to0}J_\epsilon(t)\geq J(t)$ for every $t$. Observe also that the $J_\epsilon$ are non-negative by Jensen's inequality. Hence from Fatou's Lemma it follows that 
\begin{equation*}
\liminf_{\epsilon\to0}\int_0^\infty J_\epsilon(t)dt\geq\int_0^\infty\liminf_{\epsilon\to0}J_\epsilon(t)dt\geq \int_0^\infty J(t)dt.
\end{equation*}
 This completes the proof of b).
\end{proof}
Let us indicate how, formally, the long-time asymptotics follow from the theorem. Inequality~\eqref{jensenquant} has a finite right hand side, which means that the space integrals on the left hand side have to converge to zero as $t\to\infty$. Denoting by $\n_\infty(x)$ the expected stationary state at $t=\infty$, we thus have 
\begin{equation*}
\int_0^\infty H\left(\frac{\n_\infty^a(x)}{N(x)}\right)d\mu(x)+\int_0^\infty\frac{B(x)}{N(0)}H^\infty\left(\frac{\n_\infty^s}{|\n_\infty^s|}(x)\right)|\n_\infty^s(dx)|-H\left(\int_0^\infty\frac{B(x)}{N(0)}\n_\infty(dx)\right)=0.
\end{equation*}
From the first part of the following version of Jensen's inequality it then follows that $n^a_\infty\equiv CN$ and $n^s_\infty=0$ on the support of $B$. Indeed, it suffices to set $\mu=\n(t,dx)/N(x)$, $\psi=B(x)N(x)/N(0)$ (recall that $\int B(x)N(x)dx=N(0)$) and $f=H$.
\begin{proposition}\label{jensenversion}
Let $f\in\mathcal{F}(\Rbb^n)$ be strictly convex and $\psi\geq0$ be bounded and continuous on $\overline{\Omega}$ with $\int_{\Omega}\psi(x)dx=1$.  
\begin{itemize}
\item[a)] If $\mu$ is a finite measure on $\overline{\Omega}$, then
\begin{equation}\label{jensenmeasure}
\int_\Omega \psi(x)f(\mu^a(x))dx+\int_{\overline{\Omega}}\psi(x)f^\infty\left(\frac{\mu^s}{|\mu^s|}(x)\right)d\mu^s(x)\geq f\left(\int_{\overline{\Omega}}\psi(x)d\mu(x)\right)
\end{equation}
with equality if and only if $\mu=Cdx$ on the support of $\psi$.
\item[b)] If $\{\mu_\epsilon\}_{\epsilon>0}$ is a family of finite measures on $\overline{\Omega}$ such that
\begin{equation}\label{jensenmeasurelimit}
\lim_{\epsilon\to0}\left\{\int_\Omega \psi(x)f(\mu_\epsilon^a(x))dx+\int_{\overline{\Omega}}\psi(x)f^\infty\left(\frac{\mu_\epsilon^s}{|\mu_\epsilon^s|}(x)\right)d\mu_\epsilon^s(x)- f\left(\int_{\overline{\Omega}}\psi(x)d\mu_\epsilon(x)\right)\right\}=0,
\end{equation}
then 
\begin{equation*}
\lim_{\epsilon\to0}\left(\int_\Omega|\mu_\epsilon^a(x)-m_\epsilon|\psi(x)dx+\int_{\overline{\Omega}}\psi(x)d|\mu_\epsilon^s|(x)\right)=0,
\end{equation*}
where $m_\epsilon=\int_{\overline{\Omega}}\psi(x)d\mu_\epsilon(x)$.

\end{itemize}
\end{proposition}
\begin{proof}
a) Let $\mu_\epsilon$ be a regularization of $\mu$ such that $\mu_\epsilon \stackrel{*}{\rightharpoonup}\mu$ and $\langle\mu_\epsilon\rangle\to\langle\mu\rangle$ as $\epsilon\to0$. For each $\epsilon>0$,~\eqref{jensenmeasure} follows by Jensen's inequality applied to the probability measure $\phi(x)dx$. By part b) of Proposition~\eqref{continuity}, the left hand side of~\eqref{jensenmeasure} for $\mu_\epsilon$ converges to the left hand side for $\mu$. The convergence of the integral on the right hand side as $\epsilon\to0$ follows simply from the weak* convergence of $\mu_\epsilon$ to $\mu$.

It remains to show that equality holds if and only if $\mu=Cdx$. Clearly, if $\mu=Cdx$ then equality holds. Conversely, suppose we have equality in~\eqref{jensenmeasure}. Then by the same arguments as above,
\begin{equation}\label{jensenlimit}
\lim_{\epsilon\to0}\left\{\int_\Omega \psi(x)f(\mu_\epsilon(x))dx- f\left(\int_{{\Omega}}\psi(x)\mu_\epsilon(x)dx\right)\right\}=0
\end{equation}
By an affine transformation (which may depend on $\epsilon$), we can assume w.l.o.g.\ that $\int{\mu_\epsilon\psi dx}=0$, $f(0)=0$, and $f(z)>0$ when $|z|>0$. Hence by strict convexity, for every $\delta>0$ there exists $C_\delta>0$ (depending only on $\delta$ and $f$) such that
\begin{equation*}
f(z)\geq C_\delta|z|\hspace{0.4cm}\text{whenever $|z|\geq\delta$.}
\end{equation*}
Therefore we have 
\begin{equation*}
\int_{\Omega}|\mu_\epsilon(x)|\psi(x)dx\leq\frac{1}{C_\delta}\int_\Omega f(\mu_\epsilon(x))\psi(x)dx+\delta.
\end{equation*}
This shows $\int_\Omega|\mu_\epsilon|\psi dx\to0$. If we remove our assumption $\int{\mu_\epsilon\psi dx}=0$ again, we obtain 
\begin{equation*}
\int_\Omega|\mu_\epsilon(x)-m_0|\psi(x)dx\to0,
\end{equation*} 
where $m_0:=\int_\Omega \psi(x)d\mu(x)$, and Proposition \eqref{continuity} implies
\begin{equation*}
\int_\Omega\psi(x)d|\mu(x)-m_0|=0,
\end{equation*}
whence the claim follows.

b) For each $\epsilon>0$, let $\{\mu_{\epsilon,\delta}\}_{\delta>0}$ be a regularization of $\mu_\epsilon$ such that $\mu_{\epsilon,\delta} \stackrel{*}{\rightharpoonup}\mu_\epsilon$ and $\langle\mu_{\epsilon,\delta}\rangle\to\langle\mu_\epsilon\rangle$ as $\delta\to0$. Let us denote by $J(\mu_\epsilon)$ the expression within the limit in \eqref{jensenmeasurelimit} and by $J(\mu_{\epsilon,\delta})$ the corresponding expression with $J_\epsilon$ replaced by $\mu_{\epsilon,\delta}$. Then by Proposition \eqref{continuity}, for fixed $\epsilon$ we have
\begin{equation*}
\lim_{\delta\to0}J(\mu_{\epsilon,\delta})=J(\mu_\epsilon).
\end{equation*}
Let now $\delta(\epsilon)$ be a function such that, on one hand,
\begin{equation*}
|J(\mu_{\epsilon,\delta(\epsilon)})-J(\mu_\epsilon)|<\epsilon,
\end{equation*}
and on the other hand,
\begin{equation*}
\mu_{\epsilon,\delta(\epsilon)}-\mu_\epsilon\stackrel{*}{\rightharpoonup}0\hspace{0.3cm}\text{and}\hspace{0.3cm}\langle\mu_{\epsilon,\delta(\epsilon)}-\mu_\epsilon\rangle\to0
\end{equation*}
as $\epsilon\to0$.

Then, by assumption, $\lim_{\epsilon\to0}J(\mu_{\epsilon,\delta(\epsilon)})=0$. But now we are exactly in the situation of \eqref{jensenlimit} and therefore deduce 
\begin{equation*}
\int_\Omega\psi(x)d|\mu_\epsilon(x)-m_\epsilon|\to0,
\end{equation*} 
as claimed.
\end{proof}

We are now ready to prove the following result on long-time asymptotics. We impose here in addition an assumption which guarantees that $\supp B\supset \supp\phi$, corresponding to condition (3.18) in \cite{Perthame}.

\begin{theorem}\label{alternative}
Assume in addition there exists $C>0$ such that $B(x)\geq C\phi(x)$.
Let $n^0\in\mathcal{M}^+([0;\infty))$. Then the solution of the renewal equation satisfies
\begin{equation}\label{asymp}
\lim_{t\to\infty}\int_0^\infty \phi(x)d|\n(t,x)-m_0N(x)dx|=0,
\end{equation}
where $m_0=\int_0^\infty\phi(x)dn^0(x)$.
\end{theorem}
\begin{remark} 
The assumption on $B$ can be easily relaxed to th following assumption: There exists $C>0$  and a finite set of $x_i$ ($i=1,..,N$) such that $\sum _{i=1}^N B(x-x_i)\geq C\phi(x)$. For this purpose we can observe that $H\left(\frac{\n^a(t,x)}{N(x)}\right)$ and $H^\infty\left(\frac{\n^s(t)}{|\n^s(t)|}(x)\right)|\n^s(t,dx)|$ are constant along characteristics.
\end{remark} 
\begin{proof}
Denote again by $J(t)$ the time integrand on the left hand side of \eqref{jensenquant}. As the right hand side is finite, we have $\int_0^\infty J(t)dt<\infty$, and since $J\geq0$ by Proposition \ref{jensenversion}a), we obtain a sequence $t_k\nearrow\infty$ such that $\lim_{k\to\infty}J(t_k)=0$. Apply now Proposition \ref{jensenversion}b) with $\mu_{\epsilon}$ replaced by $\n(t_k)/N$, $\psi=BN/N(0)$, and $f=H$ to obtain
\begin{equation*}
\lim_{k\to\infty}\int_0^\infty B(x)d|\n(t_k,x)-m_kN(x)|=0,
\end{equation*}
which by our assumption $B\geq C\phi$ yields
\begin{equation}\label{mk}
\lim_{k\to\infty}\int_0^\infty \phi(x)d|\n(t_k,x)-m_kN(x)|=0.
\end{equation}
Here we denoted $m_k=\int_0^\infty B(x)d\n(t_k,x)/N(0)$.

Next, let us show $\lim_{k\to\infty}m_k=m_0$. To this end, we use the conservation law
\begin{equation*}
\int_0^\infty\phi(x)d\n(t,x)=\int_0^\infty\phi(x)dn^0(x)\hspace{0.2cm}\text{for all $t>0$},
\end{equation*}
which is obtained from \eqref{contraction} by setting $H=\pm$id. But \eqref{mk} implies a fortiori
\begin{equation*}
\int_0^\infty \phi(x)d\n(t_k,x)-m_k\int_0^\infty \phi(x)N(x)\to0,
\end{equation*}
whence $m_k\to m_0$ follows since $\int\phi Ndx=1$.

Thus we have established \eqref{asymp} at least for the subsequence $\{t_k\}$. Now, consider the functional
\begin{equation*}
\mu\mapsto \int_\Omega \psi(x)f(\mu^a(x))dx+\int_{\overline{\Omega}}\psi(x)f^\infty\left(\frac{\mu^s}{|\mu^s|}(x)\right)d\mu^s(x)
\end{equation*}
for $f$ and $\psi$ as in Proposition \ref{jensenversion}. It follows from \eqref{jensenmeasure} that $m_0dx$ is the unique (up to a $\psi dx$-nullset) minimizer of this functional in the space of finite measures subject to the side constraint
\begin{equation*}
\int\psi(x)d\mu(x)=m_0.
\end{equation*}
Taking now $\psi=\phi N$, $f=H$ and $\mu=\n(t_k)/N$, and using the monotonicity property \eqref{contraction}, we find that 
\begin{equation*}
\lim_{t\to\infty}\left\{\int_0^\infty\phi(x)N(x)H\left(\frac{\n^a(t,x)}{N(x)}\right)dx+\int_0^\infty\phi(x)H^\infty\left(\frac{\n^s(t)}{|\n^s(t)|}(x)\right)|\n^s(t,dx)|\right\}=H(m_0).
\end{equation*}
Finally, still setting $\psi=\phi N$, $f=H$ and $\mu=\n(t_k)/N$, an application of Proposition \ref{jensenversion}b) yields \eqref{asymp}.\\
\end{proof}

\end{document}